\documentclass[12pt]{article}
\usepackage{amssymb,amsfonts,amsmath, psfrag,cite}
\usepackage{graphicx, colordvi}

\makeatletter \@addtoreset{equation}{section}
\@addtoreset{theo}{section} \makeatother

\def\bfn{{\boldsymbol n}}
\def\bfga{{\boldsymbol \gamma}}

\usepackage{color}

\begin{document}
\begin{center}
{\large\bf The Extended Zeilberger's Algorithm with Parameters}

\vskip 3mm

 William Y.C. Chen$^1$,
Qing-Hu Hou$^2$ and Yan-Ping Mu$^3$

\vskip 3mm
$^{1,2}$Center for Combinatorics, LPMC-TJKLC \\
Nankai University\\
Tianjin 300071, P. R. China \\

\vskip 2mm $^3$College of Science\\
Tianjin University of Technology\\
 Tianjin 300384, P.R. China \\

\vskip 2mm
 $^1$chen@nankai.edu.cn, $^2$hou@nankai.edu.cn,
$^3$yanping.mu@gmail.com

{Dedicated to Professor Wen-Tsun Wu \\ on the occasion of his
ninetieth birthday}

\end{center}

\begin{abstract}
For a hypergeometric series  $\sum_k f(k,a, b, \ldots,c)$ with
parameters $a, b, \ldots,c$, Paule has found a variation of
Zeilberger's algorithm to establish recurrence relations involving
shifts on the parameters. We consider a more general problem
concerning several similar hypergeometric terms $f_1(k, a, b,\ldots,
c)$, $f_2(k, a,b, \ldots, c)$, $\ldots$, $f_m(k, a, b, \ldots, c)$.
We present an algorithm to derive a linear relation among the sums
$\sum_k f_i(k,a,b,\ldots,c)$ $(1\leq i \leq m)$. Furthermore, when
the summand $f_i$ contains the parameter $x$, we can require that
the coefficients be $x$-free. Such relations with $x$-free
coefficients can be used to determine whether a polynomial sequence
satisfies the three term recurrence and structure relations for
orthogonal polynomials. The $q$-analogue of this approach is called
the extended $q$-Zeilberger's algorithm, which can be employed to
derive recurrence relations on the Askey-Wilson polynomials and the
$q$-Racah polynomials.
\end{abstract}

\noindent \textbf{Keywords:} Zeilberger's algorithm, the Gosper
algorithm, hypergeometric series, orthogonal polynomials

\noindent \textbf{AMS Subject Classification:}  33F10, 33C45, 33D45

\section{Introduction}

Based on Gosper's algorithm, Zeilberger \cite{Zeil91,WZ92a} has
developed a powerful theory for proving identities on hypergeometric
series and basic hypergeometric series. Let $F(n,k)$ be a double
hypergeometric term, namely, $F(n+1,k)/F(n,k)$ and $F(n,k+1)/F(n,k)$
are both rational functions of $n$ and $k$. Zeilberger's algorithm
is devised to find a double hypergeometric term $G(n,k)$ and
polynomials $a_0(n),a_1(n),\ldots,a_m(n)$ which are independent of
$k$ such that
\begin{equation}
\label{Zeil} a_0(n) F(n,k) + a_1(n) F(n+1,k) + \cdots +a_m(n)
F(n+m,k) = G(n,k+1)-G(n,k).
\end{equation}
 Writing  \[ S(n)=\sum_{k=0}^\infty F(n,k).\]
 Summing  (\ref{Zeil}) over $k$,  we deduce that
\begin{equation}\label{rec}
a_0(n) S(n) + a_1(n) S(n+1) + \cdots + a_m(n) S(n+m) =
G(n,\infty)-G(n,0).
\end{equation}
 Thus the identity
\begin{equation}
\label{iden} \sum_{k=0}^\infty F(n,k) = f(n)
\end{equation}
can be justified by verifying that $f(n)$ also satisfies \eqref{rec}
and both sides of \eqref{iden} share the same initial values.

Paule \cite{Paule05} extended Zeilberger's algorithm to the
multi-variable case and found many applications. Let $\bfn$ denote
the vector of variables $(n_1,\dots,n_r)$ and $F(\bfn,k)$ be a
multi-variable hypergeometric term, that is,
\[
\frac{F(n_1+1,n_2,\dots,n_r,k)}{F(n_1,n_2,\dots,n_r,k)}, \dots,
\frac{F(n_1,n_2,\dots,n_r+1,k)}{F(n_1,n_2,\dots,n_r,k)},
\frac{F(n_1,n_2,\dots,n_r,k+1)}{F(n_1,n_2,\dots,n_r,k)}
\]
are all rational functions of $\bfn$ and $k$. Given $m$ shifts
$\bfga_1,\dots,\bfga_m \in {\mathbb Z}^r$ of the variables $\bfn$,
he found that one can use a similar procedure to Zeilberger's
algorithm
 to find a multi-variable hypergeometric term $G(\bfn, k)$ and
coefficients $\alpha_1(\bfn), \dots, \alpha_m(\bfn)$ which are
independent of $k$ such that
\begin{equation}
\label{para} \sum_{i=1}^m \alpha_i(\bfn) F(\bfn+\bfga_i,k) =
G(\bfn,k+1)-G(\bfn,k).
\end{equation}

The main idea of this paper is the observation that Paule's approach
can be further extended to a  more general telescoping problem. Let
$f_1(k,a,b,\ldots,c)$, $\ldots$, $f_m(k,a,b,\ldots,c)$ be $m$
similar hypergeometric terms of $k$ with parameters $a,b,\ldots,c$,
namely, the ratios
\[
\frac{f_i(k,a,b,\ldots,c)}{f_j(k,a,b,\ldots,c)} \quad \mbox{and}
\quad \frac{f_i(k+1,a,b,\ldots,c)}{f_i(k,a,b,\ldots,c)}
\]
are all rational functions of $k$ and $a,b,\ldots,c$. Find a
hypergeometric term $g(k,a,b,\ldots,c)$, that is, the ratio $g(k+1,
a, b, \ldots, c)/g(k,a,b, \ldots, c)$ is a rational function of $k$
and $a,b,\ldots, c$, and polynomial coefficients
$a_1(a,b,\ldots,c)$, $a_2(a,b,\ldots,c)$, $\ldots$,
$a_m(a,b,\ldots,c)$ which are independent of $k$ such that
\begin{equation}
\label{EZ} a_1f_1(k)+a_2f_2(k)+\cdots+a_mf_m(k)=g(k+1)-g(k).
\end{equation}
For brevity,  from now on we may omit the parameters $a,b,\ldots,c$
and write $f_i(k)$ for $f_i(k,a,b,\ldots,c)$, $a_i$ for
$a_i(a,b,\ldots,c)$, and $g(k)$ for $g(k,a,b,\ldots,c)$. Once the
telescoping relation \eqref{EZ} is established, summing over $k$
 often leads to  a homogenous relation among the
sums $\sum_k f_1(k)$, $\ldots$, $\sum_k f_m(k)$:
\[
a_1 \sum_k f_1(k)+a_2 \sum_k f_2(k)+\cdots+a_m \sum_k f_m(k) = 0.
\]
Let $F(\bfn,k)$ be a multi-variable hypergeometric term and $\bfga_i
\in \mathbb{Z}^r$. Then $f_i(k,\bfn)=F(\bfn+\bfga_i,k)$ are similar
hypergeometric terms of $k$ with parameters $n_1,\ldots,n_r$.
Therefore, Paule's equation \eqref{para} is a special case of
\eqref{EZ}. However, we should note that the extended Zeilberger's
algorithm is very much in the  spirit of the original algorithm of
Zeilberger, and it should be regarded as a variation as well because
the implementation of the extended algorithm is essentially the same
as the original algorithm.

As an application of our algorithm, one can determine whether a
given hypergeometric series satisfies the recurrence relation and
the structure relations for orthogonal polynomials. Meanwhile, we
obtain the coefficients in these relations. For instance, let
$P_n(x)=\sum_k P_{n,k}(x)$ be the hypergeometric representation of
the  Jacobi polynomials as given in \eqref{jacobi}. Set
\[
f_1(k)=P_{n,k}(x), \ f_2(k)=P'_{n+1,k}(x), \ f_3(k)=P'_{n,k}(x), \
f_4(k)=P'_{n-1,k}(x),
\]
where $P_{n,k}'(x)$ denotes the derivative of $P_{n,k}(x)$ with
respect to $x$. The extended Zeilberger's algorithm enables us to
find the structure relation for $P_n(x)$
\begin{equation}
\label{abc} P_n(x) = \tilde{a}_n P'_{n+1}(x) + \tilde{b}_n P'_n(x) +
\tilde{c}_n P'_{n-1}(x).
\end{equation}
It is worth mentioning that neither the original Zeilberger's
algorithm nor the variation of Paule is directly applicable to the
above relation (\ref{abc}) involving derivatives.

Furthermore, it is important to impose an additional requirement
that the coefficients $a_1,\ldots,a_m$ in \eqref{EZ} are not only
independent of $k$ but also independent of some other parameters
such as the variable $x$ occurring as the variable of  orthogonal
polynomials. For example, $\tilde{a}_n,\tilde{b}_n$ and
$\tilde{c}_n$ in \eqref{abc} are required to be independent of the
variable $x$. Based on this parameter-free property of the
coefficients, Chen and Sun \cite{Chen-Sun} have developed a computer
algebra approach to proving identities on Bernoulli polynomials and
Euler polynomials.

We should notice  that Koepf and Schmersau \cite{Koepf-25} have
shown that one can derive the recurrence relation and structure
relations for orthogonal polynomials by variations of Zeilberger's
algorithm. For each of the three kinds of relations,  they have
provided an algorithm.  The extended Zeilberger's algorithm serves
as a unification of their algorithms and applies to more general
cases.  For instance, our algorithm can also be used to derive
recurrence relations for the connection coefficients between two
classes of Meixner polynomials with different parameters.

In another direction, the extended Zeilberger's algorithm can be
adapted to deal with basic hypergeometric terms. Using the
$q$-analogue of this algorithm, we can recover the three term
recurrence relations for the Askey-Wilson polynomials and the
$q$-Racah polynomials.

Let us recall some terminology and notation. A function $t(k)$ is
called a hypergeometric term if $t(k+1)/t(k)$ is a rational function
of $k$. A hypergeometric series is defined by
\[
{_rF_s}  \left( \left. \begin{array}{c} a_1,\ldots,a_r\\
b_1,\ldots,b_s
\end{array} \right| z \right) = \sum_{k=0}^\infty \frac{(a_1)_k
\cdots (a_r)_k}{(b_1)_k  \cdots (b_s)_k} \frac{z^k}{k!},
\]
where $(a)_k=a(a+1) \cdots (a+k-1)$ is the raising factorial. The
$q$-shifted factorial is given by
$(a;q)_k=(1-a)(1-aq)\cdots(1-aq^{k-1})$ and we write
\[
(a_1,\ldots,a_m;q)_k = (a_1;q)_k \cdots (a_m;q)_k.
\]
Then a basic hypergeometric series is defined by
\[
{_r\phi_s} \left[ \left. \begin{array}{c} a_1,\ldots,a_r\\
b_1,\ldots,b_s
\end{array} \right| q; z \right] =
\sum_{k=0}^\infty \frac{(a_1, \cdots, a_r;q)_k} {(b_1, \cdots
b_s;q)_k} {z^k \over (q;q)_k} \left( (-1)^k q^{k \choose 2}
\right)^{s-r+1}.
\]

\section{The Extended Zeilberger's Algorithm}

Let $f_1(k), f_2(k), \ldots, f_m(k)$ be similar hypergeometric terms
with parameters $a, b, \ldots, c$. Recall that two hypergeometric
terms $f(k)$ and $g(k)$ are said to be similar if their ratio is a
rational function of $k$ and the parameters. We assume that
\begin{equation} \label{fij}
\frac{f_1(k+1)}{f_1(k)}= \frac{u(k)}{v(k)} \quad \mbox{and} \quad
\frac{f_i(k)}{f_1(k)} = \frac{p_i(k)}{Q(k)}, \ i=1,2,\ldots,m,
\end{equation}
where $u(k),v(k),p_i(k),Q(k)$ are polynomials in $k$ and the
parameters $a, b, \ldots, c$. Suppose that $f_i(k)$ satisfy
\eqref{fij}. Then
\[
\frac{f_i(k+1)}{f_i(k)} = \frac{f_i(k+1)/f_1(k+1)}{f_i(k)/f_1(k)}
\frac{f_1(k+1)}{f_1(k)} =\frac{p_i(k+1) Q(k)u(k)}{p_i(k) Q(k+1)v(k)}
\]
and
\[
\frac{f_i(k)}{f_j(k)} = \frac{f_i(k)/f_1(k)}{f_j(k)/f_1(k)} =
\frac{p_i(k)}{p_j(k)}
\]
are rational functions of $k$ and $a, b,\ldots, c$. Thus
\eqref{fij} is equivalent to the statement that $f_1(k), f_2(k),
\ldots, f_m(k)$ are similar hypergeometric terms.

Our aim is to find coefficients $a_1,\ldots,a_m$ as rational
functions in the parameters $a, b, \ldots, c$ but independent of $k$
(called {\it $k$-free} coefficients) such that
\begin{equation}
\label{EZ-a} a_1f_1(k)+a_2f_2(k)+\cdots+a_mf_m(k) = g(k+1)-g(k)
\end{equation}
for some hypergeometric term $g(k)$ with parameters $a, b, \ldots,
c$. By the similarity of $f_1(k),\ldots,f_m(k)$,
\begin{equation}\label{tk}
t_k = a_1f_1(k)+a_2f_2(k)+\cdots+a_mf_m(k)
\end{equation}
is a hypergeometric term of $k$ with parameters $a, b, \ldots,c$. So
we can apply Gosper's algorithm \cite{Gosper78}  to find $g(k)$ such
that $t_k=g(k+1)-g(k)$. Notice that we always have a trivial
solution $a_1=a_2=\cdots=a_m=0$ and $g(k)=0$.

Notice that by multiplying the common denominator, the coefficients
$a_1,\ldots,a_m$ in \eqref{EZ-a} become polynomials in the
parameters $a,b,\ldots,c$.  If no confusion arises, we may not
mention the parameters $a, b, \ldots, c$. In the usual case $n$ is
the parameter for identities on finite sums.

It follows from (\ref{tk}) that
\begin{eqnarray*}
\frac{t_{k+1}}{t_k} &=& \frac{f_1(k+1)}{f_1(k)} \frac{ \sum_{i=1}^m
a_i f_i(k+1)/f_1(k+1)} {\sum_{i=1}^m a_i f_i(k)/f_1(k)} \\
&=& \frac{u(k) Q(k)}{v(k)Q(k+1)} \frac{\sum_{i=1}^m a_i
p_i(k+1)}{\sum_{i=1}^m a_i p_i(k)}.
\end{eqnarray*}
Suppose that
\[
\frac{u(k) Q(k)}{v(k)Q(k+1)} = \frac{a(k)}{b(k)} \frac{c(k+1)}{c(k)}
\]
is a Gosper representation, i.e., $a(k),b(k),c(k)$ are polynomials
such that $\gcd(a(k),b(k+h))=1$ for all non-negative integers $h$.
Then a Gosper representation of $t_{k+1}/t_k$ is given by
\[
\frac{t_{k+1}}{t_k} = \frac{a(k)}{b(k)}
\frac{c(k+1)P(k+1)}{c(k)P(k)},
\]
where
\begin{equation}\label{Pk}
P(k) = \sum_{i=1}^m a_i p_i(k).
\end{equation}
Gosper's algorithm states that $g(k)$ exists if and only if there
exists a polynomial $x(k)$ such that
\begin{equation}
\label{G-key} a(k)x(k+1) - b(k-1) x(k) = c(k)  P(k) .
\end{equation}
Moreover, the degree bound $d$ for $x(k)$ can be estimated by $a(k)$
and $b(k)$. Suppose \[ x(k)=\sum_{i=0}^d c_i k^i.\] By comparing the
coefficients of each power of $k$ on both sides, we obtain a system
of linear equations on $a_1,\ldots,a_m$ and $c_0,c_1,\ldots, c_d$.
Solving the system of linear equations, we find coefficients
$a_1,\ldots,a_m$ and
\[
g(k)= \frac{b(k-1)x(k)}{c(k)Q(k)} f_1(k).
\]

The extended Zeilberger's algorithm can be described in terms of the
following steps.

\noindent {\it Input:} $m$ similar hypergeometric terms $f_1(k),\ldots,f_m(k)$. \\
{\it Output:} $k$-free coefficients $a_1,a_2,\ldots,a_m$ and a
hypergeometric term $g(k)$ satisfying  \eqref{EZ-a}.

\begin{itemize}
\item[1.]
Compute the rational functions \[ r_i(k) =
\frac{f_i(k)}{f_1(k)}\quad\mbox{and}\quad
r(k)=\frac{f_1(k+1)}{f_1(k)}.\] Set $Q(k)$ to be the common
denominator of $r_1(k),\ldots,r_m(k)$, \[ p_i(k)=r_i(k) Q(k),\] and
let $P(k)$ be given by \eqref{Pk}.
\item[2.]
Compute a Gosper representation of
\[
r(k)\frac{Q(k)}{Q(k+1)} = \frac{a(k)}{b(k)} \frac{c(k+1)}{c(k)}.
\]
\item[3.]
Compute the degree bound $d$ for $x(k)$ and solve the equation
\eqref{G-key} by the method of undetermined coefficients to obtain
the $k$-free coefficients $a_1,\ldots,a_m$ and the polynomial
$x(k)$.
\item[4.]
The hypergeometric term $g(k)$ is then given by
\[
g(k)= \frac{b(k-1)x(k)}{c(k)Q(k)} f_1(k).
\]
\end{itemize}

Suppose that $F(n,k)$ is a double hypergeometric term. Let $f_i(k) =
F(n+i-1,k)$. Then the extended Zeilberger's algorithm reduces to the
original Zeilberger's algorithm. More generally, suppose that
$\bfn=(n_1,\ldots,n_r)$, $F(\bfn,k)$ is a multi-variable
hypergeometric term and $\bfga_i \in \mathbb{Z}^r$. The
specialization $f_i(k)=F(\bfn+\bfga_i,k)$ reduces to Paule's
variation.

As will be seen, in some applications it is necessary to require
that the coefficients $a_1,\ldots,a_m$ be independent of some
parameters, say, the parameter $a$. More precisely, let
$f_1(k,a,b,\ldots,c)$, $\ldots$, $f_m(k,a,b,\ldots,c)$ be $m$
similar multi-variable hypergeometric terms, that is,
\begin{equation}
\label{para-rat} \frac{f_i(k+1,a,b,\ldots,c)}{f_i(k,a,b,\ldots,c)}
\quad \mbox{and} \quad
\frac{f_i(k,a,b,\ldots,c)}{f_1(k,a,b,\ldots,c)}
\end{equation}
are rational functions of $k,a,b,\ldots,c$. We aim to find $a_1,a_2,
\ldots,a_m$ not only independent of $k$ but also independent of the
parameter $a$ such that \eqref{EZ-a} holds.

Since the solutions $(a_1,\ldots,a_m,g(k))$ of \eqref{EZ-a} form a
linear vector space, we may use the following form as the output
\begin{align}
& a_1 = v_1, \ldots, a_r=v_r, \nonumber \\
& a_{r+1}=h_{r+1}(v_1,\ldots,v_r), \ldots, a_m =
 h_m(v_1,\ldots,v_r), \label{hh}\\
& g(k)=h(v_1,\ldots,v_r) f_1(k), \nonumber
\end{align}
where $v_1,\ldots,v_r$ are free variables and $h_{r+1},\ldots,h_m,
h$ are linear combinations of $v_1,\ldots,v_r$. For this purpose, we
should first ignore the independence of $a$ and apply the extended
Zeilberger's algorithm to find the solution \eqref{hh}. By
\eqref{para-rat}, the functions $h_{r+1},\ldots,h_m$ are rational
functions of the parameters $a, b, \ldots, c$ and thus can be
written as
\[
h_i=p_i(a,b, \ldots, v_1,\ldots,v_r)/q_i(a, b, \ldots, c),
\]
where $p_i,q_i$ are relatively prime polynomials. Now consider the
additional requirement that $a_1,a_2,\ldots,a_m$ are independent of
the parameters $a$. Therefore,  all the coefficients of
\[
p_i(a,b, \ldots, c, v_1,\ldots,v_r) - a_i q_i(a, b, \ldots, c)
\]
in variable $a$ must be zero. This gives a system of linear
equations on $a_1,\ldots,a_m$ and $v_1,\ldots,v_r$. Upon solving
these equations, we eventually find $a_1,a_2,\ldots,a_m$ which are
independent of $k$ and the parameter $a$. The above version of the
extended Zeilberger's algorithm will still be called the extended
Zeilberger's algorithm.

Let us take the Hermite polynomials as the first example to show how
to use the above algorithm to derive linear relations on sums of
similar hypergeometric terms with parameters.

\noindent {\bf Example 2.1} The Hermite polynomials $H_n(x)$ are
given by
\begin{equation}\label{Hn}
 H_n(x) = (2x)^n \ {_2F_0} \left( \left.
\begin{array}{c} -\frac{n}{2}, -\frac{n-1}{2} \\ \rule{10pt}{0.5pt}
\end{array} \right| -\frac{1}{x^2} \right),
\end{equation}
see \cite[Section 6.1]{AAR99}. We aim to find a  three term
recurrence
 \[ xH_n(x) = \alpha_n H_{n+1}(x) + \beta_n H_n(x) +
\gamma_n H_{n-1}(x)
\]
with  coefficients $\alpha_n,\beta_n,\gamma_n$ being independent of
$x$.  Let
\[
H_{n,k}(x) = (2x)^n \frac{\big( -\frac{n}{2} \big)_k \big(
\frac{-n+1}{2} \big)_k}{k!} \left(-\frac{1}{x^2} \right)^k
\]
be the summand in \eqref{Hn}. We first ignore the $x$-freeness
requirement and apply the extended Zeilberger's algorithm to the
four similar hypergeometric terms with parameters $n$ and $x$
\[
f_1(k)=xH_{n,k}(x), \ f_2(k)=H_{n+1,k}(x), \ f_3(k)=H_{n,k}(x), \
f_4(k)=H_{n-1,k}(x).
\]
We find that
\begin{equation}
\label{ex1} a_1=v_1, \ a_2=v_2, \ a_3=-x(v_1+2v_2), \ a_4=2n v_2,
\end{equation}
and
\begin{equation} \label{ex1b}
g(k) = \frac{-4k v_2}{n+1-2k} x H_{n,k}(x).
\end{equation}
Now it is time to impose  the $x$-freeness condition to give  an
additional equation
\[ v_1+2v_2=0.\] Hence we obtain
\[
a_1=v_1, \ a_2=-\frac{v_1}{2}, \ a_3=0, \ a_4=-n v_1, \ g(k) =
\frac{2k v_1}{n+1-2k} xH_{n,k}(x).
\]
It follows that
\[
v_1xH_{n,k}(x) - \frac{v_1}{2} H_{n+1,k}(x) - n v_1 H_{n-1,k}(x) =
g(k+1)-g(k).
\]
Summing over $k$, we deduce that
\[
xH_n(x) = \frac{1}{2} H_{n+1}(x) + n H_{n-1}(x).
\]

\section{Orthogonal Polynomials}

Using the extended Zeilberger's algorithm, we can determine whether
a hypergeometric series satisfies a three term relation and the
structure relations for a sequence of orthogonal polynomials. In
other words, one can verify the orthogonality of a terminating
hypergeometric series by using the extended Zeilberger's algorithm.

The method to derive the relation in Example 2.1 is in fact valid in
the general case. Given a hypergeometric series $P_n(x)$, we can
compute the coefficients for the following recurrence relation
\[
xP_n(x) = \alpha_n P_{n+1}(x) + \beta_n P_n(x) + \gamma_n
P_{n-1}(x).
\]
 Let
 $ P_n(x)=\sum_k
P_{n,k}(x)$, where $P_{n,k}(x)$ is a hypergeometric term of $k$ with
parameters $n$ and $x$. Set
\[
f_1(k)=xP_{n,k}(x),\ f_2(k)=P_{n+1,k}(x),\ f_3(k)=P_{n,k}(x),\
f_4(k)=P_{n-1,k}(x).
\]
Clearly, $f_1, f_2, f_3$ and $f_4$ are similar hypergeometric terms.
So we can use the extended Zeilberger's algorithm to find $a_1, a_2,
a_3,a_4$ which are independent of $k$ and $x$ such that
\[
a_1 xP_{n,k}(x) + a_2 P_{n+1,k}(x) + a_3 P_{n,k}(x) + a_4
P_{n-1,k}(x) = g(k+1)-g(k).
\]
Then summing over $k$ leads to
\[
\alpha_n = -a_2/a_1,\ \beta_n = -a_3/a_1,\ \gamma_n = -a_4/a_1.
\]

By this method, we can recover the three term recurrences for the
Laguerre polynomials, the Jacobi polynomials, the Charlier
polynomials, the Meixner polynomials, the Kravchuk polynomials and
the Hahn polynomials, as listed in the following table. Notice that
we have adopted the notation in \cite{KS98}.
 \allowdisplaybreaks
\begin{center}
\begin{tabular}{c|l}
monic OPs & coefficients $(\alpha_n=1)$ \\
\hline \rule{0pt}{30pt}

\begin{tabular}{c}
Laguerre \\
$L^{(a)}_n(x)$
\end{tabular}
&
$\beta_n=a+2n+1, \  \gamma_n=n(a+n)$ \\[15pt]

\begin{tabular}{c}
Jacobi \\
$P^{(a,\,b)}_n(x)$
\end{tabular} & $ \begin{cases}
{ \beta_n=\frac{-(a-b)(a+b)}{(2n+2+a+b)(2n+a+b)}}
\\[8pt]
  { \gamma_n =
\frac{4(n+b)(a+n)(n+a+b)n}{(2n+a+b+1)(2n+a+b-1)(2n+a+b)^2}}
\end{cases}$ \\[25pt]

\begin{tabular}{c}
Charlier \\
$C_n(x;a)$
\end{tabular} & $\beta_n=a+n, \ \gamma_n = an$ \\[20pt]

\begin{tabular}{c}
Meixner \\
$M_n(x;b,c)$
\end{tabular} & $\displaystyle \beta_n=\frac{cb+nc+n}{1-c}, \ \gamma_n =
 \frac{nc(b+n-1)}{(c-1)^2}$ \\[20pt]

\begin{tabular}{c}
Krawtchouk \\
$K_n(x;p,N)$
\end{tabular} & $\begin{cases} \beta_n=Np-2np+n, \\ \gamma_n=pn(1-p)(N-n+1) \end{cases}$ \\[25pt]

\begin{tabular}{c}
Hahn \\
$Q_n(x;a,b,N)$
\end{tabular} & $ \begin{cases} \beta_n=\frac{
(b+a+2n^2+2n+2nb+2na+a^2+ab)N-n(a-b)(n+a+b+1)}{(2n+a+b)(2n+2+a+b)}
\\[8pt]
 \gamma_n =
\frac{n(N-n+1)(n+b)(a+n)(n+a+b)(a+b+N+1+n)}{(2n+a+b+1)(2n+a+b-1)(2n+a+b)^2}
\end{cases}
$
\end{tabular}
\end{center}

\noindent {\bf Example 3.1} The Wilson polynomials $W_n(x)$
 are given by
\[
\frac{W_n(x^2)}{(a+b)_n (a+c)_n (a+d)_n} = {_4F_3} {\small \left(
\left.
\begin{array}{c} -n,n+a+b+c+d-1,a+xi,a-xi \\ a+b,a+c,a+d \end{array} \right| 1
\right)},
\]
where $i=\sqrt{-1}$, see \cite{Wilson}. Let $W_{n,k}(x^2)$ be the
summand of the right hand side multiplied by the factor $(a+b)_n
(a+c)_n (a+d)_n$. Applying the the extended Zeilberger's algorithm
to the similar terms
\[
 f_1(k)=x^2W_{n,k}(x^2),\ f_2(k)=W_{n+1,k}(x^2),\ f_3(k)=W_{n,k}(x^2),\ f_4(k)=W_{n-1,k}(x^2),
\]
we obtain the following relation, see also \cite[p.\,24]{KS98}
\[
x W_n(x) = \alpha_n W_{n+1}(x) + \beta_n W_n(x) + \gamma_n
W_{n-1}(x),
\]
where
\[
\alpha_n= - \frac{a+b+c+d+n-1}{(a+b+c+d+2n)(a+b+c+d+2n-1)},
\]
\begin{multline*}
\beta_n  = \frac{4n^2-4(1-a-b-c-d)n+(a+b+c+d)^2-2(a^2+b^2+c^2+d^2)}{8} \\[10pt]
  + \frac{ (a+b+c+d-2)(a+b-c-d)(a+c-b-d)(a+d-b-c)}
{8(a+b+c+d+2n)(a+b+c+d+2n-2)},
\end{multline*}
and
\begin{multline*}
\gamma_n = - (a+b+n-1)(a+c+n-1)(a+d+n-1) \\[5pt]
\times \frac{
(b+c+n-1)(b+d+n-1)(c+d+n-1)n}{(a+b+c+d+2n-1)(a+b+c+d+2n-2)}.
\end{multline*}

\noindent {\bf Example 3.2} The Racah polynomials $R_n(x)$
 are given by
\[
R_n(x(x+c+d+1)) = {_4F_3} {\small \left( \left.
\begin{array}{c} -n,n+a+b+1,-x,x+c+d+1 \\ a+1,b+d+1,c+1 \end{array} \right| 1
\right)},
\]see \cite{A-W79}.
The extended Zeilberger's algorithm gives the recurrence relation of
the Racah polynomials \cite[p.\,27]{KS98}
\[
x R_n(x) = \alpha_n R_{n+1}(x) + \beta_n R_n(x) + \gamma_n
R_{n-1}(x),
\]
where
\[
\alpha_n =
\frac{(a+n+1)(c+n+1)(a+b+n+1)(d+b+n+1)}{(a+b+2n+1)(a+b+2n+2)},
\]
\begin{multline*}
\beta_n = \frac{-4n^2-4(a+b+1)n +(a-2d-b-2-2c)(a-b)}{8}\\[5pt]
 -\frac{(c+1)(b+d+1)}{2}
 -\frac{(a-b)(a+b)(a-2d-b)(a-2c+b)}{8(a+b+2n)(a+b+2n+2)},
\end{multline*}
and
\[
\gamma_n = \frac{ (b+n)(a-d+n)(a+b-c+n)n}{(a+b+2n)(a+b+2n+1)}.
\]

\noindent {\bf Example 3.3} Askey and Ismail \cite{Askey84} provided
two hypergeometric representations for the Pollaczek polynomials:
\begin{eqnarray*}
P_n(x) &=& \eta^n \, {_2F_1} \left( \left.
\begin{array}{c} -n,b(x-\eta)/\xi, \\ b/a  \end{array} \right|
-\frac{\xi}{a \eta} \right) \\
&=& \zeta^n \, {_2F_1} \left( \left.
\begin{array}{c} -n,b(\zeta-x)/\xi, \\ b/a  \end{array} \right|
\frac{\xi}{a \zeta} \right),
\end{eqnarray*}
where
\[
\xi = \sqrt{(1+a)^2x^2 - 4a}, \quad \eta = ((1+a)x - \xi)/2a, \quad
\zeta=((1+a)x + \xi)/2a.
\]
Using the extended Zeilberger's algorithm, we derive the following
three term recurrence from either  representation \cite{Askey84}:
\[
xP_n(x) = \frac{an+b}{(1+a)n+b} P_{n+1}(x) + \frac{n}{(1+a)n+b}
P_{n-1}(x).
\]

We continue to show that  the extended Zeilberger's algorithm can be
employed to express the derivatives of  orthogonal polynomials in
terms the original polynomials, and vice versa. Let $P_{n,k}(x)$ be
the summand of the hypergeometric representation of $P_n(x)$ and
$P'_{n,k}(x)$ be the derivative of $P_{n,k}(x)$. It is easily seen
that $P'_{n,k}(x)$ is similar to $P_{n,k}(x)$. This enables us to
derive the three term recurrence for $P'_n(x)$ and the structure
relations for $P_n(x)$ as given below
\begin{equation}\label{s-1}
 \sigma(x) P'_n(x) =
a_n P_{n+1}(x) + b_n P_{n}(x) + c_n P_{n-1}(x),
\end{equation}
and
\begin{equation}\label{s-2}
P_n(x) =  \bar{a}_n P'_{n+1}(x) + \bar{b}_n P'_{n}(x) + \bar{c}_n
P'_{n-1}(x),
\end{equation}
where $\sigma(x)$ is a polynomials in $x$ of degree less than or
equal to $2$ and $a_n,b_n,c_n,\bar{a}_n,\bar{b}_n,\bar{c}_n$ are
constants not depending on $x$. To derive \eqref{s-1}, we set
\[
f_1(k)=\sigma(x) P'_{n,k}(x), f_2(k)=P_{n+1,k}(x),\
f_3(k)=P_{n,k}(x),\ f_4(k)=P_{n-1,k}(x).
\]
To establish \eqref{s-2}, we set
\[
f_1(k)=P_{n,k}(x), f_2(k)=P'_{n+1,k}(x),\ f_3(k)=P'_{n,k}(x),\
f_4(k)=P'_{n-1,k}(x).
\]

\noindent {\bf Example 3.4} The monic Jacobi polynomials
 are given by
\begin{equation}\label{jacobi}
P_n(x) = \frac{(a+1)_n2^n}{(n+a+b+1)_n} \ {_2F_1} \left( \left.
\begin{array}{c} -n,n+a+b+1 \\ a+1 \end{array} \right| \frac{1-x}{2}
\right).
\end{equation}
Let $P_{n,k}(x)$ denote the summand. Its derivative with respect to
$x$ equals
\[
P'_{n,k}(x)= - \frac{(a+1)_n2^n}{(n+a+b+1)_n} \frac{(-n)_k
(n+a+b+1)_k}{2(a+1)_k (k-1)!} \left(\frac{1-x}{2} \right)^{k-1}.
\]
Consider the four similar terms
\[
f_1(k)=x P'_{n,k}(x), \ f_2(k)=P'_{n+1,k}(x), \ f_3(k)=P'_{n,k}(x),
\ f_4(k)=P'_{n-1,k}(x).
\]
By the extended Zeilberger's algorithm with parameters $n$ and $x$,
we find that
\begin{multline*}
x P'_n(x) = \frac{n}{n+1} P'_{n+1}(x) -
\frac{(a+2+b)(a-b)}{(2n+2+a+b)(2n+a+b)} P'_n(x) \\[5pt]
+ \frac{4n(b+n)(a+n)(n+a+b+1)}{(2n+a+b+1)(2n+a+b-1)(2n+a+b)^2}
P'_{n-1}(x),
\end{multline*}
\begin{multline*}
(1-x^2)P'_n(x) = -n P_{n+1}(x)
+\frac{2n(a-b)(n+a+b+1)}{(2n+2+a+b)(2n+a+b)} P_n(x) \\[5pt]
+\frac{4(n+b)(a+n)(n+a+b+1)(n+a+b)n}{(2n+a+b+1)(2n+a+b-1)(2n+a+b)^2}
P_{n-1}(x).
\end{multline*}
and
\begin{multline*}
P_n(x) = \frac{1}{n+1} P'_{n+1}(x) +
\frac{2(a-b)}{(2n+2+a+b)(2n+a+b)} P'_n(x) \\[5pt]
\qquad \quad - \frac{4(n+b)(a+n)n}{(2n+a+b+1)(2n+a+b-1)(2n+a+b)^2}
P'_{n-1}(x).
\end{multline*}

The following example is concerned with expressing the shifts of
orthogonal polynomials with parameters in terms of the original
polynomials and their derivatives.

\noindent {\bf Example 3.5} Let
\[
P_n^{(a,\,b)}(x) = \frac{(a+1)_n2^n}{(n+a+b+1)_n} \ {_2F_1} \left(
\left.
\begin{array}{c} -n,n+a+b+1 \\ a+1 \end{array} \right| \frac{1-x}{2}
\right)
\]
be the Jacobi polynomials, see \cite{AAR99, KS98}. By applying the
extended Zeilberger's algorithm to $f_1(k)=P_{n,k}^{(a+1,\,b)}(x)$
($f_1(k)=P_{n,k}^{(a,\,b+1)}(x)$, respectively) and
\[
f_2(k)={P_{n+1,k}^{(a,\,b)}}\,'(x),\
f_3(k)={P_{n,k}^{(a,\,b)}}\,'(x),\
f_4(k)={P_{n-1,k}^{(a,\,b)}}\,'(x),
\]
we are led to the known relations due to Koepf and Schmersau
\cite{Koepf-25}
\[
P_n^{(a+1,\,b)}(x) = \frac{1}{n+1} {P_{n+1}^{(a,\,b)}}\,'(x) +
\frac{2(a+1+n)}{(2n+2+a+b)(2n+a+b+1)} {P_n^{(a,\,b)}}\,'(x)
\]
and
\[
P_n^{(a,\,b+1)}(x) = \frac{1}{n+1} {P_{n+1}^{(a,\,b)}}\,'(x) -
\frac{2(b+1+n)}{(2n+2+a+b)(2n+a+b+1)} {P_n^{(a,\,b)}}\,'(x).
\]
Moreover, we can  deduce the following relations which seem to be
new:
\begin{multline*} P_n^{(a+1,\,b-1)}(x) = \frac{1}{n+1}
{P_{n+1}^{(a,\,b)}}\,'(x) +
\frac{4(a+1+n)}{(2n+2+a+b)(2n+a+b)} {P_{n}^{(a,\,b)}}\,'(x) \\[5pt]
\hfill + \frac{4(a+1+n)(a+n)n}{(2n+a+b-1)(2n+a+b+1)(2n+a+b)^2}
{P_{n-1}^{(a,\,b)}}\,'(x),
\end{multline*}
and
\begin{multline*} P_n^{(a-1,\,b+1)}(x) = \frac{1}{n+1}
{P_{n+1}^{(a,\,b)}}\,'(x) -
\frac{4(b+1+n)}{(2n+2+a+b)(2n+a+b)} {P_{n}^{(a,\,b)}}\,'(x) \\[5pt]
\hfill + \frac{4(b+1+n)(b+n)n}{(2n+a+b-1)(2n+a+b+1)(2n+a+b)^2}
{P_{n-1}^{(a,\,b)}}\,'(x).
\end{multline*}

The extended Zeilberger's algorithm can also be employed to compute
the connection coefficients of two sequences of orthogonal
polynomials. Ronveaux \cite{Ronveaux} developed an approach to
computing recurrence relations for the connection coefficients. The
extended Zeilberger's algorithm serves this purpose  without
resorting to the properties of the orthogonal polynomials. As an
example, let us consider the connection coefficients of two classes
of Meixner polynomials with different parameters.

\noindent {\bf Example 3.6} Let $M_n^{(\gamma,\,\mu)}(x)$ be the
monic Meixner polynomials defined by
\[
M_n^{(\gamma,\,\mu)}(x) = (\gamma)_n \left( \frac{\mu}{\mu-1}
\right)^n \ {_2F_1} \left( \left.
\begin{array}{c} -n,-x \\ \gamma
\end{array} \right| 1-\frac{1}{\mu} \right),
\] see \cite[p.\,45]{KS98}.
We wish to find a recurrence relation for the connection
coefficients $C_m(n)$ defined by
\begin{equation}\label{M-M}
M_n^{(\gamma,\,\mu)}(x) = \sum_{m=0}^n C_m(n)
M_m^{(\delta,\,\nu)}(x).
\end{equation}
To this end, we first find a difference operator which eliminates
$M_n^{(\gamma,\,\mu)}(x)$. This goal can be achieved by applying the
extended Zeilberger's algorithm to the similar terms
\[
f_1(k)=M_{n,k}^{(\gamma,\,\mu)}(x), \
f_2(k)=M_{n,k}^{(\gamma,\,\mu)}(x+1), \ \mbox{and}\
f_3(k)=M_{n,k}^{(\gamma,\,\mu)}(x-1),
\]
where
\[
M_{n,k}^{(\gamma,\,\mu)}(x) = (\gamma)_n \left( \frac{\mu}{\mu-1}
\right)^n \frac{(-n)_k (-x)_k}{(\gamma)_k k!} \left( 1-\frac{1}{\mu}
\right)^k.
\]
From the telescoping relation generated by the extended Zeilberger's
algorithm, we deduce that
\[
(x\mu+\mu\gamma+x-n+n\mu) M_n^{(\gamma,\,\mu)}(x) -\mu (\gamma+x)
M_n^{(\gamma,\,\mu)}(x+1) -x M_n^{(\gamma,\,\mu)}(x-1) = 0.
\]
Let
\[
S_m(x) = (x\mu+\mu\gamma+x-n+n\mu) M_m^{(\delta,\,\nu)}(x) -\mu
(\gamma+x) M_m^{(\delta,\,\nu)}(x+1) -x M_m^{(\delta,\,\nu)}(x-1),
\]
which can be used to establish a linear relation on the connection
coefficients $C_m(n)$. Indeed, it follows from \eqref{M-M}  that
\begin{equation}\label{cm}
\sum_{m=0}^n C_m(n) S_m(x) = 0.
\end{equation}
Suppose that we can express $S_m(x)$ in terms of a suitable basis
$\{P_m(x)\}$:
\begin{equation}\label{S-P}
S_m(x) = a_m P_{m+1}(x) + b_m P_m(x) + c_m P_{m-1}(x),
\end{equation}
where $a_m,b_m$ and $c_m$ are independent of $x$. Substituting
(\ref{S-P}) into (\ref{cm}), by the linear independence of $P_m(x)$
for $m=0, 1, 2, \ldots$, that is, the coefficients of $P_i(x)$ are
all zeros, we find
\begin{equation}\label{C-rec}
a_{m-1} C_{m-1}(n) + b_m C_m(n) + c_{m+1} C_{m+1}(n) = 0.
\end{equation}
Thus the question has become how to find the polynomials $P_m(x)$ in
order to determine the coefficients $a_m$, $b_m$ and $c_m$. In view
of the relation \eqref{S-P}, we consider a hypergeometric term
$P_m(x)$ that is  similar to $S_m(x)$ so that we can solve the
equation
\[ S_m(x)-a_mP_{m+1}(x)   - b_m P_m(x) -c_m P_{m-1}(x)=0\] by using
the extended Zeilberger's algorithm. In fact, we may choose
\[
P_m(x)=\Delta(M_m^{(\delta,\,\nu)}(x)) = M_m^{(\delta,\,\nu)}(x+1) -
M_m^{(\delta,\,\nu)}(x).
\]
It is easily checked that $P_m(x)$ satisfies \eqref{S-P} and the
corresponding coefficients are given by
\[
a_m =\frac{(\mu-1)(n-m)}{m+1}, \quad c_m=\frac{(\nu-\mu)(\delta+m-1)
m \nu}{(1-\nu)^2},
\]
\[
b_m = \frac{-\nu\mu m - m\mu + 2m\nu + \nu\mu\gamma-\nu n + \nu
\delta - \nu + \mu - \nu\mu\delta - \mu \gamma + \nu n \mu}{1-\nu}.
\]
Hence we have derived a recurrence relation  \eqref{C-rec} for the
connection coefficients $C_m(n)$.

\section{$q$-Orthogonal Polynomials}

The extended Zeilberger's algorithm can be readily adapted to basic
hypergeometric terms $t_k$ with parameters $a, b, \ldots, c$, that
is, the ratio of two consecutive terms is a rational function of
$q^k$ and the parameters. The $q$-analogue of the extended
Zeilberger's algorithm will be called the extended $q$-Zeilberger's
algorithm. Let $f_1(k)$, $f_2(k), \ldots, f_m(k)$ be similar
$q$-hypergeometric terms, namely,
\[
f_i(k)/f_j(k)  \quad \mbox{and} \quad f_i(k+1)/f_i(k)\ (1 \le i, j
\le m)
\]
are rational functions of $q^k$ and the parameters. The objective of
the extended $q$-Zeilberger's algorithm is  to find a
$q$-hypergeometric term $g(k)$ and coefficients $a_1$,
$a_2,\ldots,a_m$ which are independent of $k$ such that
\begin{equation}
\label{qEZ} a_1 f_1(k)+a_2 f_2(k)+\cdots+a_m f_m(k)=g(k+1)-g(k).
\end{equation}
The detailed description of the extended $q$-Zeilberger's algorithm
is similar to that of the ordinary case, hence is omitted. We will
only give examples to demonstrate how to use this method to compute
the three term recurrences and structure relations for
$q$-orthogonal polynomials.

\noindent {\bf Example 5.1} The discrete $q$-Hermite I polynomials
are given by \cite{AC65}
\[
H_n(x) = q^{n \choose 2} \ {_2\phi_1} \left[ \left.
\begin{array}{c} q^{-n}, x^{-1} \\[6pt]  0 \end{array} \right| q;-qx
\right].
\]
Let $dH_n(x) = \frac{H(xq)-H(x)}{(q-1)x}$ be the $q$-difference of
$H_n(x)$. We derive that
\[
x dH_n(x) = \frac{1-q^n}{1-q^{n+1}} dH_{n+1}(x) + q^{n-2}(1-q^n)
dH_{n-1}(x)
\]
and
\[
dH_n(x) = \frac{1-q^n}{1-q}H_{n-1}(x).
\]

\noindent {\bf Example 5.2} The Askey-Wilson polynomials
 $p_n(x;a,b,c,d|q)$ are defined  by
\[
\frac{a^n p_n(x;a,b,c,d|q)}{(ab,ac,ad;q)_n} = \ {_4\phi_3} \left[
\left.
\begin{array}{c}
q^{-n},abcdq^{n-1},ae^{i\theta},ae^{-i\theta}
\\[6pt] ab,ac,ad \end{array} \right| q;q
\right], \ x= \cos \theta,
\]
 see \cite[(7.5.2)]{GasRah04}.
Let $t_{n,k}(x)$ be the summand of the right hand side multiplied by
$(ab,ac,ad;q)_n/a^n$. Applying the  extended $q$-Zeilberger's
algorithm to
\[
f_1(k)=xt_{n,k}(x),\ f_2(k)=t_{n+1,k}(x),\ f_3(k)=t_{n,k}(x),\
f_4(k)=t_{n-1,k}(x),
\]
we find that
\[
x p_n(x) = \alpha_n p_{n+1}(x) + \beta_n p_n(x) + \gamma_n
p_{n-1}(x),
\]
where
\[
\alpha_n =  \frac{1-abcdq^{n-1}}{2 (1-abcdq^{2n})(1-abcdq^{2n-1})},
\]
\begin{multline*} \beta_n =
\frac{q^{n-1} (abcd q^{2n-1}+1) ((a+b+c+d)q+bcd+acd+abd+abc)}{2(1-abcdq^{2n})(1-abcd q^{2n-2})}
\\[5pt]
   -\frac{q^{2n-2} (1+q)((bcd+acd+abd+abc)q+abcd(a+b+c+d))
}{2(1-abcdq^{2n})(1-abcd q^{2n-2})},
\end{multline*}
and
\begin{multline*}
\gamma_n= (1-q^n)(1-abq^{n-1})(1-acq^{n-1})(1-adq^{n-1}) \\[5pt]
\times \frac{(1-bcq^{n-1})(1-bdq^{n-1})(1-cdq^{n-1})}
  {2(1-abcdq^{2n-1})(1-abcdq^{2n-2})}.
\end{multline*}

\noindent {\bf Example 5.3} The $q$-Racah polynomials
$R_n(x;a,b,c,d|q)$ are given by
\[
R_n(q^{-x}+cdq^{x+1};a,b,c,d|q) = \ {_4\phi_3} \left[ \left.
\begin{array}{c}
q^{-n},abq^{n+1},q^{-x},cdq^{x+1}
\\[6pt] aq,bdq,cq \end{array} \right| q;q
\right],
\]
see \cite[p.\,122]{KS98}. The extended $q$-Zeilberger's algorithm
gives the following recurrence relation first derived by Askey and
Wilson \cite{A-W79} using a transformation formula on a ${}_8\phi_7$
series:
\[ x R_n(x) = \alpha_n
R_{n+1}(x) + \beta_n R_n(x) + \gamma_n R_{n-1}(x),
\]
where
\[
\alpha_n =
\frac{(1-aq^{n+1})(1-abq^{n+1})(1-bdq^{n+1})(1-cq^{n+1})}{(1-abq^{2n+1})(1-abq^{2n+2})},
\]
\begin{multline*}
\beta_n = \frac{q^{n+1}(abq^{2n+1}+1) (c+bcd+dc+a+bd+ab+ca+abd)}
{(1-ab q^{2n})(1-ab q^{2n+2})} \\[5pt]
- \frac{q^{2n+1} (1+q) (ab^2d+abcd+ca+bcd+abd+ab+abc+a^2b)}{(1-ab
q^{2n})(1-ab q^{2n+2})},
\end{multline*}
and
\[
\gamma_n=
\frac{(1-q^n)(1-bq^n)(c-abq^n)(d-aq^n)q}{(1-abq^{2n})(1-abq^{2n+1})}.
\]

\vskip 15pt \noindent {\small {\bf Acknowledgments.}  This work was
supported by the 973 Project on Mathematical Mechanization, the
National Natural Science Foundation, the PCSIRT project of the
Ministry of Education, and the Ministry of Science and Technology of
China. Y.-P. Mu was supported by National Natural Science Foundation
of China, Project 10826038.}

\end{document}